\documentclass[final]{siamltex}

\usepackage{epsfig}
\usepackage{fancyhdr}
\usepackage{multirow}
\usepackage{amsmath,amssymb}
\usepackage {graphicx}
\usepackage{algorithm}
\usepackage{algorithmic}
\newtheorem{thm}{Theorem}[section]

\title{Semi-Randomized Block Kaczmarz Methods with Simple Random Sampling for Large-Scale Linear Systems}

\author{Gang Wu\thanks{Corresponding author. 
School of Mathematics, China University of Mining and Technology, 221116, Jiangsu, P.R. China. E-mail: {\tt gangwu@cumt.edu.cn}. }
       \and Qiao Chang\thanks{ School of Mathematics,
China University of Mining and Technology, Xuzhou, 221116, Jiangsu, P.R. China.
E-mail: {\tt 1597064283@qq.com}.}
}

\begin{document}

\maketitle

\begin{abstract}
Randomized block Kaczmraz method plays an important role in solving large-scale linear system. One of the key points of this type of methods is how to effectively select working rows. However, in most of the state-of-the-art randomized block Kaczmarz-type methods, one has to scan all the rows of the coefficient matrix in advance for computing probabilities or paving, or to compute the residual vector of the linear system in each iteration to determine the working rows. Thus, we have to access all the rows of the data matrix in these methods, which are unfavorable for big-data problems. Moreover, to the best of our knowledge, how to efficiently choose working rows in randomized block Kaczmarz-type methods for multiple linear systems is still an open problem. In order to deal with these problems, we propose semi-randomized block Kaczmarz methods with simple random sampling for linear systems with single and multiple right-hand sides, respectively. In these methods, there is no need to scan or pave all the rows of the coefficient matrix, nor to compute probabilities and the residual vector of the linear system in each step. Specifically, one can update all the solutions of a large-scale linear system with multiple right-hand sides simultaneously. The convergence of the proposed methods is considered. Numerical experiments on both real-world and synthetic data sets show that the proposed methods are superior to many state-of-the-art randomized Kaczmarz-type methods for large-scale linear systems.
\end{abstract}



\section{Introduction}
Consider the large-scale consistent linear system with {\it multiple} right-hand sides
\begin{equation}\label{1.1}
 \mathbf{AX}=\mathbf{B},
\end{equation}
where $\mathbf{A}\in R^{m\times n}$, $\mathbf{X}\in R^{n\times k_{b}}$, $\mathbf{B}\in R^{m\times k_{b}}$ and $k_{b}\geq 1$. Specifically, if $k_{b}=1$, then \eqref{1.1} reduces
to the large-scale linear system with {\it simple} right-hand side
\begin{equation}\label{1.2}
 \mathbf{Ax}=\mathbf{b}.
\end{equation}

The Kaczmarz method is popular for solving large-scale consistent linear systems, and it is known as an algebraic reconstruction technique. This method has been widely used in medical image scanning \cite{9}, tomography \cite{1,10,15,23}, and so on. The original Kaczmarz method was proposed in \cite{13}. In this method, one selects the working rows sequentially, and then project the current solution to the hyperplane solution space determined by the row, to generate the next approximate solution
\begin{equation*}
\mathbf{x}_{t+1}=\mathbf{x}_{t}+\frac{\mathbf{\mathbf{b}}_{(i)}-\mathbf{A}_{(i)}\mathbf{x}_{t}}{\|\mathbf{A}_{(i)}\|_{2}^{2}}\mathbf{A}_{(i)}^{T},\quad t=0,1,\ldots,
\end{equation*}
where~$\mathbf{b}_{(i)}$ and $\mathbf{A}_{(i)}$ denote the $i$-th row of ~$\mathbf{b}$ and $\mathbf{A}$, respectively, and $t$ is the number of iterations.

In \cite{24}, Strohmer and Vershyn proposed a randomized Kaczmarz (RK) method, which is a pioneered research on this topic. The key is to choose the working rows in terms of the probability generated by the proportion of row norms. The convergence with expected exponential rate was also established.
The RK method is quite appealing in practical applications, and has received great
attention from many researchers. For example, Hefny {\it et al.} presented variants of
the RK method and the randomized Gauss-Siedel (RGS) method \cite{HNR} for solving
large linear systems. A weighted randomized Kaczmarz method was presented in \cite{St}, and a randomized extended Kaczmarz method for solving the least squares problems were proposed in \cite{25}. In \cite{3}, Bai and Wu proposed an efficient greedy randomized Kaczmarz method that improves the randomized Kaczmarz method significantly. However, it is required to calculate the residual of the linear system and construct an indicator set in each iteration, which is unfavorable for large-scale data problems. In order to partially deal with this problem, Jiang {\it et al.} \cite{12} proposed a semi-randomized Kaczmarz method with simple random sampling. 



As a generalization of the randomized Kaczmarz method, the block Kaczmarz methods have also been investigated extensively in recent years \cite{4,5,8,16,17,18,19,20,21}. For instance, Needell and Tropp \cite{20} proposed a randomized block Kaczmarz (RBK) method. Niu and Zheng \cite{21} considered a greedy block Kaczmarz method by using a greedy strategy. 
In this method, there is no need to determine the partition of the row index of the coefficient matrix in advance.
In \cite{16}, Liu and Gu  proposed a greedy randomized block Kaczmarz method (GRBK) and a block extension of Motzkin method (BEM). 
Necoara \cite{18} proposed a randomized average block Kaczmarz (RABK) method. In \cite{5}, Du {\it et al.} proposed a simple stochastic extended average block Kaczmarz method (REABK), which converges exponentially to the unique minimum norm least-squares solution of a given linear system. An advantage of the average block Kaczmarz-type methods is they are pseudo-inverse free, however, it is often required to choose many parameters in this type of methods. 

All the randomized (block) Kaczmarz methods mentioned above are just for the large linear system with simple right-hand side \eqref{1.2}. To the best of our knowledge, how to efficiently choose working rows in randomized block Kaczmarz-type
methods for multiple linear systems is still an open problem. Moreover, in most of the above methods, one has to scan all the rows of the coefficient matrix $\mathbf{A}$ for paving \cite{20}, compute the residual vector to construct indicator sets during iterations, or determine some parameters in advance. So it is interesting to give new insight into the randomized block Kaczmarz method.

In this work, we pay special attention to block Kaczmarz method for large-scale linear system with simple and multiple right-hand sides. 
However, generalizing randomized (block) Kaczmarz method from simple right-hand side problem to multiple right-hand sides problems is not straightforward. The key lies in how to achieve multiple right-hand solutions at the same time, and select working rows effectively. 
Motivated by the semi-randomized Kaczmarz method with simple random sampling \cite{12}, we propose a partially randomized block Kaczmarz method for large-scale linear systems with simple and multiple right-hand sides. Specifically, one can update all the solutions of large-scale linear system with multiple right-hand sides simultaneously.  
Moreover, we do not need to scan all the rows of the data matrix for paving, nor to compute all the elements of the residual vector during iterations. The convergence of the proposed methods is given.

The rest of this paper is organized as follows. In Section 2, we present some preliminaries of this work. In Sections 3, we propose our
semi-randomized block Kaczmarz methods with simple random sampling for simple right-hand side and multiple right-hand sides, respectively. The convergence results are given. Numerical experiments are performed in Section 4, and some concluding remarks are given in Section 5.

Given the matrix $\mathbf{A}\in R^{ m\times n}$, let $\mathbf{A}^{T}$, $\mathbf{A}^{\dagger}$, and $\|\mathbf{A}\|_{F}$ be the transpose, the Moore-Penrose inverse, the Frobenius norm of ~$\mathbf{A}$, respectively. Let $\mathbf{b}_{(i)}$ be the $i$th component of the vector $\mathbf{b}$, and $\mathbf{A}_{(i)}$ be the $i$-th row of $\mathbf{A}$. For any random variable $\mathbf{x}$, we denote by $\mathbb{E}[\mathbf{x}]$ its expectation, and by $Nu(\mathbf{A})$ the number of rows of $\mathbf{A}$. In this paper, $\mathbf{A}_{Jt}$
stands for the row submatrix of $\mathbf{A}$ indexed by $J_{t}$, and $Nu(\mathbf{A}_{Jt})$ represents the number of rows in $\mathbf{A}_{Jt}$. Let $\mathbf{B}_{(:,j)}$ be the $j$-th column of $\mathbf{B}$, and $\mathbf{B}_{(i,j)}$ be the element of $\mathbf{B}$ in the $i$-th row and the $j$-th column.


 \section{Two randomized Kaczmarz-type methods for large-scale linear systems}
In this section, we present some preliminaries, and briefly introduce  
the randomized block Kaczmarz method \cite{20} as well as the semi-randomized Kaczmarz method with simple random sampling \cite{12}.

 \subsection{Randomized block Kaczmarz method}
Given a partition $\mathcal{T}=\{\tau_{1},\tau_{2},\ldots,\tau_{p}\}$ of the row indices $\{1,2,\ldots,m\}$, the coefficient matrix $\mathbf{A}$ and the right-hand side $\mathbf{b}$ can be divided into $p$ blocks as follows
 \begin{equation}
\mathbf{A}^{T}
=
\begin{pmatrix}
\mathbf{A}_{\tau_{1}}^{T},\mathbf{ A}_{\tau_{2}}^{T}, \ldots, \mathbf{A}_{\tau_{p}}^{T}
\end{pmatrix},
\end{equation}
and
\begin{equation}
\mathbf{b}^{T}
=
\begin{pmatrix}
\mathbf{b}_{\tau_{1}}^{T}, \mathbf{b}_{\tau_{2}}^{T}, \ldots, \mathbf{b}_{\tau_{p}}^{T}
\end{pmatrix},
\end{equation}
where $\mathbf{A}_{\tau_{i}}\in R^{m_{i}\times n}$, $\mathbf{b}_{\tau_{i}}\in R^{m_{i}}$ are the sub-matrix and sub-vector of $\mathbf{A}$ and $\mathbf{b}$ with components listed in $\tau_{i}$, respectively, and $\sum\limits_{i=1}\limits^{p} m_{i}=m$.
Different from the non-block Kaczmarz-type methods, the block Kaczmarz method is equivalent to solving multiple equations in each iteration \cite{8}. More precisely, if we choose $\mathbf{A}_{\tau_{i}}$ as the working rows, then the current iterative solution~$\mathbf{x}_{t}$ is projected onto the hyperplane~$\mathbf{A}_{\tau_{i}}\mathbf{x}=\mathbf{b}_{\tau_{i}}$. In other words,
$$
\mathbf{x}_{t+1}=\mathbf{x}_{t}+\mathbf{A}_{\tau_{i}}^{\dag}(\mathbf{b}_{\tau_{i}}-\mathbf{A}_{\tau_{i}}\mathbf{x}_{t}).
$$

In randomized block Kaczmarz methods \cite{5,16,17,18,19,20,21},
the choice of a block at each iteration is crucial in practice. For instance, one can choose a block uniformly at random from $\mathcal{T}$,
or select $\tau_{i}$ with probability $pr(\tau_{i})=\frac{\|\mathbf{A}_{\tau_{i}}\|_{F}^{2}}{\|\mathbf{A}\|_{F}^{2}}$. A framework of the randomized block Kaczmarz algorithm (RBK) is described as follows.
 \begin{algorithm}\label{Alg1}
\caption{A randomized block Kaczmarz method (RBK) \cite{20}.}


${\rm \bf Input}$:

\ \  $\mathbf{A},\ \mathbf{b},\ l$, $\mathbf{x}_{0}$, and partition $\mathcal{T}=\{\tau_{1},\tau_{2},\ldots,\tau_{p}\}$ of the row indices of $\{1,2,\ldots,m\}$.

{\rm \bf Output}: Approximate solution $\hat{\mathbf{x}}$.

1. For $t=0,1,\ldots,l-1$

2.~~~~ Select a block $\tau_{i}$ from $\mathcal{T}$.

3.~~~~ Update $\mathbf{x}_{t+1}=\mathbf{x}_{t}+\mathbf{A}_{\tau_{i}}^{\dag}(\mathbf{b}_{\tau_{i}}-\mathbf{A}_{\tau_{i}}\mathbf{x}_{t})$.

4.~~~~ If $\mathbf{x}_{t+1}$ is accurate enough, then let $\hat{\mathbf{x}}=\mathbf{x}_{t+1}$ and stop, else continue.

5. End for
\end{algorithm}

\subsection{A semi-randomized Kaczmarz method with simple random sampling}

Recently, a semi-randomized Kaczmarz method with simple random sampling was proposed in \cite{12}, for solving the large-scale linear system with simple right-hand side \eqref{1.2}.
Unlike some existing Kaczmarz methods, there is no need to calculate probabilities and residuals of the linear system, nor to construct index sets via scanning the residuals in this method. Indeed, one only needs to compute some elements of residual vectors corresponding to simple sampling sets, and thus a small portion of rows of the coefficient matrix are utilized.



More precisely, the key idea is that one first chooses a small set $\Omega_t$ in the $t$-th iteration by using a simple random sampling method, and then
seeks the working row in $\Omega_t$ rather than in all the $m$ rows of $A$.
In this method, one chooses the working row $A_{i_t}$ using the following probability
\begin{equation}\label{23}
\frac{|\mathbf{r}_{(i_t)}|}{\|\mathbf{A}_{(i_t)}\|_{2}}=\frac{|\mathbf{b}_{(i_t)}-\mathbf{A}_{(i_t)}\mathbf{x}_{t}|}{\|\mathbf{A}_{(i_t)}\|_{2}}=\max\limits_{j \in \Omega_{t}}\left\{\frac{|\mathbf{b}_{(j_t)}-\mathbf{A}_{(j_t)}\mathbf{x}_{t}|}{\| \mathbf{A}_{(j_t)}\|_{2}}\right\},
\end{equation}
where $\mathbf{r}_{(i_t)}$~is the~$i$-th element of~ the residual $\mathbf{r}_t=\mathbf{b}-\mathbf{Ax}_t$~of the~$t$-th iteration.
The semi-randomized Kaczmarz method with simple random sampling is given below.



\begin{algorithm}\label{Alg2}
\caption{A semi-randomized Kaczmarz method with simple random sampling for large-scale linear system with {\it simple} right-hand side (SRK) \cite{12}.}

${\rm \bf Input}$:

\ \  $\mathbf{A},\ \mathbf{b}, 0<\eta\ll 1, \ l$, and $\mathbf{x}_{0}$.

{\rm \bf Output}: Approximate solution $\hat{\mathbf{x}}$.

1. For $t=0,1,\ldots,l-1$ do

~~~~2. Generating an indicator set $\Omega_{t}$, i.e., choosing $\eta m$ rows of $A$ by using the simple random sampling method.

~~~~3. Select the working row according to \eqref{23}.

~~~~4. Let $\mathbf{x}_{t+1}=\mathbf{x}_{t}+\frac{\mathbf{b}_{(i_t)}-\mathbf{A}_{(i_t)}\mathbf{x}_{t}}{\|\mathbf{A}_{(i_t)}\|_{2}^{2}}\mathbf{A}_{(i_t)}^{T}$.

~~~~5. If $\mathbf{x}_{t+1}$ is accurate enough, then let $\hat{\mathbf{x}}=\mathbf{x}_{t+1}$ and stop, else continue.

6. End for
\end{algorithm}



\section{Semi-randomized block Kaczmarz methods with simple random sampling}

In this section, we propose semi-randomized block Kaczmarz methods for large-scale linear systems with {\it simple} and {\it multiple} right-hand sides, respectively. Different from some existing randomized block Kaczmarz methods \cite{5,8,16,17,18,19,20,21}, we construct an indicator set rather than choose it from $\mathcal{T}$ in each iteration.

\subsection{Solving large-scale linear system with {\it simple} right-hand side}


In this subsection, we focus on the large-scale linear system \eqref{1.2} with {\it simple} right-hand side. In randomized block Kaczmarz methods, one of the key points is to choose working rows in each iteration. At the $t$-th iteration, we denote by
\begin{equation}\label{3.1}
\beta_{(j_t)}=\frac{|\mathbf{b}_{(j_t)}-\mathbf{A}_{(j_t)}\mathbf{x}_{t}|}{\|\mathbf{A}_{(j_t)}\|_{2}},\quad j\in\{1,2,\ldots,m\}.
\end{equation}
Motivated by \eqref{23}, the idea is to select a row index set $J_{t}$ from $A$ corresponding to some large values of the $\{\beta_{(j_t)}\}$'s. That is, we construct 
a row index set $J_t$ from $A$ corresponding to the largest $k$ values of the $\{\beta_{(j_t)}\}$'s. More precisely, suppose that 
\begin{equation}\label{31}
\beta_{({j^{(1)}_t})}\geq\beta_{({j^{(2)}_t})}\geq\cdots\geq\beta_{({j^{(k)}_t})}
\end{equation}
are the $k$ largest ones in  $\{\beta_{(j_t)}\}$'s, then we set
\begin{equation}\label{32}
J_t=\{j^{(1)},j^{(2)},\ldots,j^{(k)}\}.
\end{equation}
Here $k\ll m$ is the number of working rows or the size of the indicator set, which can be chosen artificially. The semi-randomized {\it block} Kaczmarz method for
large-scale linear system with {\it simple} right-hand side is given in Algorithm 3.


\begin{algorithm}
\caption{A semi-randomized block Kaczmarz method for
large-scale linear system with {\it simple} right-hand side}


${\rm \bf Input}$:

\ \  $\mathbf{A},\ \mathbf{b},\ l$,  $\mathbf{x}_{0}$, and $k_{max}$.

{\rm \bf Output}:

\ \  $\mathbf{x}_{l}$ .

1. For $t=0,1,\ldots,l-1$ do

~~~~2. Compute $\beta_{(jt)}=\frac{|\mathbf{b}_{(jt)}-\mathbf{A}_{(jt)}\mathbf{x}_{t}|}{\|\mathbf{A}_{(jt)}\|_{2}}$ with $jt\in\{1,2,\ldots,m\}$.

~~~~3. Sort the $\{\beta_{(jt)}\}$'s and select the $k_{max}$ largest ones. Construct the index set $J_{t}$ according to \eqref{32}. 

~~~~4. Update~$\mathbf{x}_{t+1}=\mathbf{x}_{t}+\mathbf{A}_{J_{t}}^{\dag}(\mathbf{b}_{J_{t}}-\mathbf{A}_{J_{t}}\mathbf{x}_{t})$.

~~~~5. If $\mathbf{x}_{t+1}$ is accurate enough, then stop, else continue.

6. End for
\end{algorithm}

We are ready to consider the convergence of the above algorithm.
\begin{thm}\label{Thm3.1}
The iteration sequence $\{\mathbf{x}_{t}\}_{t=0}^{\infty}$  generated by  Algorithm 3 converges to the least-squares solution $\mathbf{x}_{\star}=\mathbf{A}^{\dagger}\mathbf{b}$, with
\begin{equation}
\|\mathbf{x}_{t+1}-\mathbf{x}_{\star}\|_{2}^{2}\leq\left( 1-\varsigma_{t}\right)\|\mathbf{x}_{t}-\mathbf{x}_{\star}\|_{2}^{2}
\end{equation}
where
\begin{equation}
 \varsigma_{t}=\frac{\|\mathbf{A}_{J_{t}}\|_{F}^{2}}{\|\mathbf{A}\|_{F}^{2}+ (\lambda_{t}-1)\|\mathbf{A}_{J_{t}}\|_{F}^{2}}\frac{ \sigma_{\min}^{2}(\mathbf{A}) }{\sigma_{\max}^{2}(\mathbf{A}_{J_{t}})},
\end{equation}
and
\begin{equation}\label{35}
\lambda_{t}=\max\limits_{jt\in J_{t}}\left\{\frac{\frac{|\mathbf{b}_{(jt)}-\mathbf{A}_{(jt)}\mathbf{x}_{t}|^{2}}{\|\mathbf{A}_{(jt)}\|_{2}^{2}}}{{\min\limits_{tt\in J_{t}}\left\{\frac{|\mathbf{b}_{(tt)}-\mathbf{A}_{(tt)}\mathbf{x}_{t}|^{2}}{\|\mathbf{A}_{(tt)}\|_{2}^{2}}\right\}}} \right\}\geq 1,\quad jt\in \{1,2,\ldots,m\}.
\end{equation}
Here $\sigma_{\min}(\mathbf{A})$, $\sigma_{\max}(\mathbf{A}_{J_{t}})$ denote the smallest nonzero and the largest nonzero
singular value of $A$ and $\mathbf{A}_{J_{t}}$, respectively.
\end{thm}
\begin{proof}
We have from Algorithm 3 that
$$
\mathbf{x}_{t+1}=\mathbf{x}_{t}+\mathbf{A}_{J_{t}}^{\dag}(\mathbf{b}_{J_{t}}-\mathbf{A}_{J_{t}}\mathbf{x}_{t}),\quad t=0,1,\ldots
$$
Let $\mathbf{x}_{\star}=\mathbf{A}^{\dagger}\mathbf{b}$ be the least-squares solution of \eqref{1.2}, then
$$
\mathbf{x}_{t+1}-\mathbf{x}_{\star}=\mathbf{x}_{t}-\mathbf{x}_{\star}+\mathbf{A}_{J_{t}}^{\dag}(\mathbf{b}_{J_{t}}-\mathbf{A}_{J_{t}}\mathbf{x}_{t}).
$$
From the fact that $\mathbf{A}_{J_{t}}^{\dag}\mathbf{A}_{J_{t}}$ is an orthogonal projector, we obtain
$$
\|\mathbf{x}_{t}-\mathbf{x}_{\star}\|_{2}^{2}=\|\mathbf{A}_{J_{t}}^{\dag}\mathbf{A}_{J_{t}}(\mathbf{x}_{t}-\mathbf{x}_{\star})\|_{2}^{2}+\|\mathbf{x}_{t+1}-\mathbf{x}_{\star}\|_{2}^{2}.
$$
That is,
\begin{equation}\label{36}
\|\mathbf{x}_{t+1}-\mathbf{x}_{\star}\|_{2}^{2}=\|\mathbf{x}_{t}-\mathbf{x}_{\star}\|_{2}^{2}-\|\mathbf{A}_{J_{t}}^{\dag}\mathbf{A}_{J_{t}}(\mathbf{x}_{t}-\mathbf{x}_{\star})\|_{2}^{2}.
\end{equation}

It follows that
\begin{equation}\label{388}
\|\mathbf{A}_{J_{t}}^{\dag}\mathbf{A}_{J_{t}}(\mathbf{x}_{t}-\mathbf{x}_{\star})\|_{2}^{2}\geq\sigma_{\min}^{2}(\mathbf{A}_{J_{t}}^{^{\dag}})\sum\limits_{i=1}\limits^{k_{max}}|\mathbf{b}_{(it)}-\mathbf{A}_{(it)}\mathbf{x}_{t}|^{2}.
\end{equation}
where we used 
$$
\frac{|\mathbf{b}_{(it)}-\mathbf{A}_{(it)}\mathbf{x}_{t}|^{2}}{\|\mathbf{A}_{(it)}\|_{2}^{2}}\geq\min\limits_{tt\in J_{t}}\left\{\frac{|\mathbf{b}_{(tt)}-\mathbf{A}_{(tt)}\mathbf{x}_{t}|^{2}}{\|\mathbf{A}_{(tt)}\|_{2}^{2}}\right\},\quad i=1,2,\ldots,k_{max}.
$$
Therefore,
\begin{small}\begin{align}\label{39}
   \|\mathbf{A}_{J_{t}}^{\dag}\mathbf{A}_{J_{t}}(\mathbf{x}_{t}-\mathbf{x}_{\star})\|_{2}^{2}
   \geq&\sigma_{\min}^{2}(\mathbf{A}_{J_{t}}^{^{\dag}})\min\limits_{tt\in J_{t}}\left\{\frac{|\mathbf{b}_{(tt)}-\mathbf{A}_{(tt)}\mathbf{x}_{t}|^{2}}{\|\mathbf{A}_{(tt)}\|_{2}^{2}}\right\}\sum\limits_{it\in J_{t}}\|\mathbf{A}_{(it)}\|_{2}^{2}\\ \nonumber
   =&\sigma_{\min}^{2}(\mathbf{A}_{J_{t}}^{^{\dag}})\|\mathbf{A}_{J_{t}}\|_{F}^{2}\min\limits_{tt\in J_{t}}\left\{\frac{|\mathbf{b}_{(tt)}-\mathbf{A}_{(tt)}\mathbf{x}_{t}|^{2}}{\|\mathbf{A}_{(tt)}\|_{2}^{2}}\right\}\\ \nonumber
   =&\frac{\|\mathbf{A}_{J_{t}}\|_{F}^{2}}{\sigma_{\max}^{2}(\mathbf{A}_{J_{t}})}\min\limits_{tt\in J_{t}}\left\{\frac{|\mathbf{b}_{(tt)}-\mathbf{A}_{(tt)}\mathbf{x}_{t}|^{2}}{\|\mathbf{A}_{(tt)}\|_{2}^{2}}\right\}\frac{\| \mathbf{A}\mathbf{x}_{\star}-\mathbf{A}\mathbf{x}_{t}\|_{2}^{2}}{\| \mathbf{b}-\mathbf{A}\mathbf{x}_{t}\|_{2}^{2}}\\ \nonumber
   \geq&\frac{\|\mathbf{A}_{J_{t}}\|_{F}^{2}}{\sigma_{\max}^{2}(\mathbf{A}_{J_{t}})}\sigma_{\min}^{2}(\mathbf{A})\frac{\min\limits_{tt\in J_{t}}\left\{\frac{|\mathbf{b}_{(tt)}-\mathbf{A}_{(tt)}\mathbf{x}_{t}|^{2}}{\|\mathbf{A}_{(tt)}\|_{2}^{2}}\right\} \|\mathbf{x}_{t}-\mathbf{x}_{\star}\|_{2}^{2}}{\left(\sum\limits_{jt\in J_{t}}\frac{|\mathbf{b}_{(jt)}-\mathbf{A}_{(jt)}\mathbf{x}_{t}|^{2}}{\|\mathbf{A}_{(jt)}\|_{2}^{2}} \|\mathbf{A}_{(jt)}\|_{2}^{2} +\sum\limits_{jt\in [m]/J_{t}}\frac{|\mathbf{b}_{(jt)}-\mathbf{A}_{(jt)}\mathbf{x}_{t}|^{2}}{\|\mathbf{A}_{(jt)}\|_{2}^{2}} \|\mathbf{A}_{(jt)}\|_{2}^{2}\right)} \\ \nonumber
   \geq&\frac{\sigma_{\min}^{2}(\mathbf{A})}{\sigma_{\max}^{2}(\mathbf{A}_{J_{t}})}\|\mathbf{A}_{J_{t}}\|_{F}^{2}\frac{\|\mathbf{x}_{t}-\mathbf{x}_{\star}\|_{2}^{2}}{\sum\limits_{jt\in J_{t}}\frac{\frac{|\mathbf{b}_{(jt)}-\mathbf{A}_{(jt)}\mathbf{x}_{t}|^{2}}{\|\mathbf{A}_{(jt)}\|_{2}^{2}}}{\min\limits_{tt\in J_{t}}\left\{\frac{|\mathbf{b}_{(tt)}-\mathbf{A}_{(tt)}\mathbf{x}_{t}|^{2}}{\|\mathbf{A}_{(tt)}\|_{2}^{2}}\right\}}\|\mathbf{A}_{(jt)}\|_{2}^{2}+\sum\limits_{jt\in [m]/  J_{t}}\frac{\frac{|\mathbf{b}_{(jt)}-\mathbf{A}_{(jt)}\mathbf{x}_{t}|^{2}}{\|\mathbf{A}_{(jt)}\|_{2}^{2}}}{\min\limits_{tt\in J_{t}}\left\{\frac{|\mathbf{b}_{(tt)}-\mathbf{A}_{(tt)}\mathbf{x}_{t}|^{2}}{\|\mathbf{A}_{(tt)}\|_{2}^{2}}\right\}}\|\mathbf{A}_{(jt)}\|_{2}^{2}}.
   \end{align}
\end{small}
Here $\mathbf{A}_{[m]/ [J_{t}]}$ is the submatrix obtained from removing $\mathbf{A}_{J_{t}}$ from $\mathbf{A}$, and the row index set corresponding to $\mathbf{A}_{[m]/ [J_{t}]}$ is denoted by $[m]/ [J_{t}]$.

On the other hand, as
$$\frac{\frac{|\mathbf{b}_{(jt)}-\mathbf{A}_{(jt)}\mathbf{x}_{t}|^{2}}{\|\mathbf{A}_{(jt)}\|_{2}^{2}}}{\min\limits_{tt\in J_{t}}\left\{\frac{|\mathbf{b}_{(tt)}-\mathbf{A}_{(tt)}\mathbf{x}_{t}|^{2}}{\|\mathbf{A}_{(tt)}\|_{2}^{2}}\right\}}\leq 1,\quad jt\in[m]/ [J_{t}],$$
we get 
\begin{small}\begin{align}\nonumber
   \|\mathbf{A}_{J_{t}}^{\dag}\mathbf{A}_{J_{t}}(\mathbf{x}_{t}-\mathbf{x}_{\star})\|_{2}^{2}
   \geq&\frac{\sigma_{\min}^{2}(\mathbf{A})}{\sigma_{\max}^{2}(\mathbf{A}_{J_{t}})}\|\mathbf{A}_{J_{t}}\|_{F}^{2}\frac{\|\mathbf{x}_{t}-\mathbf{x}_{\star}\|_{2}^{2}}{\sum\limits_{jt\in J_{t}}\lambda_{t}\|\mathbf{A}_{(jt)}\|_{2}^{2}+\sum\limits_{jt\in[m]/J_{t}}\|\mathbf{A}_{(jt)}\|_{2}^{2}} \\\nonumber
   \geq&\frac{\sigma_{\min}^{2}(\mathbf{A})}{\sigma_{\max}^{2}(\mathbf{A}_{J_{t}})}\|\mathbf{A}_{J_{t}}\|_{F}^{2}\frac{\|\mathbf{x}_{t}-\mathbf{x}_{\star}\|_{2}^{2}}{ (\lambda_{t}-1)\|\mathbf{A}_{J_{t}}\|_{F}^{2}+\|\mathbf{A}\|_{F}^{2}}.
   \end{align}
\end{small}

Moreover, we notice that
\begin{equation}\label{38}
\|\mathbf{A}_{J_{t}}^{\dag}\mathbf{A}_{J_{t}}(\mathbf{x}_{t}-\mathbf{x}_{\star})\|_{2}^{2}
\geq\frac{\|\mathbf{A}_{J_{t}}\|_{F}^{2}}{\|\mathbf{A}\|_{F}^{2}+ (\lambda_{t}-1)\|\mathbf{A}_{J_{t}}\|_{F}^{2}}\frac{ \sigma_{\min}^{2}(\mathbf{A}) }{\sigma_{\max}^{2}(\mathbf{A}_{J_{t}})}\|\mathbf{x}_{t}-\mathbf{x}_{\star}\|_{2}^{2}.
\end{equation}
Substituting \eqref{38} into \eqref{36}, we arrive at
$$
\|\mathbf{x}_{t+1}-\mathbf{x}_{\star}\|_{2}^{2}\leq\left( 1-\frac{\|\mathbf{A}_{J_{t}}\|_{F}^{2}}{\|\mathbf{A}\|_{F}^{2}+ (\lambda_{t}-1)\|\mathbf{A}_{J_{t}}\|_{F}^{2}}\frac{ \sigma_{\min}^{2}(\mathbf{A}) }{\sigma_{\max}^{2}(\mathbf{A}_{J_{t}})}\right)\|\mathbf{x}_{t}-\mathbf{x}_{\star}\|_{2}^{2},
$$
which completes the proof.
\end{proof}

However, one needs to compute the $\{\beta_{(jt)}\}$'s, $jt=1,2,\ldots,m$; see \eqref{3.1}.
That is, we have to compute all the elements of the residual vector at each iteration. In other words, it is required to access the matrix $\mathbf{A}$. This goes against the motivation of the original Kaczmarz method.

Motivated by the idea advocated in \cite{12}, we randomly select a subset $\mathbf{A}_{t},\mathbf{b}_{t}$ from $\mathbf{A}$ and $\mathbf{b}$ by using the technique of simple random sampling. Let us discuss it in more detail. Let $\Omega_{t}$ be a subset of $\{1,2,\ldots,m\}$, which is generated by using simple random sampling, with $Nu(\Omega_{t})=\eta m$, $0<\eta\ll 1$. Thus, instead of \eqref{3.1}, we compute 
\begin{equation}
\beta_{(jt)}=\frac{|\mathbf{b}_{(jt)}-\mathbf{A}_{(jt)}\mathbf{x}_{t}|}{\|\mathbf{A}_{(jt)}\|_{2}},~~jt\in \Omega_{t}.
\end{equation}
That is, we only need to calculate $\eta m(\ll m)$ elements of the residual vector, rather than the $m$ ones. 
By using this strategy, 
there is no need to calculate all the elements in the residual vector, nor scan all the data matrix $\mathbf{A}$.
Compared with Algorithm 3, one can reduce the workload and save the storage requirements significantly, especially for big data problems. We are ready to present the following semi-randomized {\it block} Kaczmarz method with simple random sampling for large-scale linear system with {\it simple} right-hand side.
\begin{algorithm}\label{Alg4}
\caption{A semi-randomized {\it block} Kaczmarz method with simple random sampling for
large-scale linear system with {\it simple} right-hand side}


${\rm \bf Input}$:

\ \  $\mathbf{A},\ \mathbf{b},\ l$,  $\mathbf{x}_{0}$, $k_{max}$, $0<\eta\ll 1$.

{\rm \bf Output}:

\ \  $\mathbf{x}_{l}$ .

1. For $t=0,1,\ldots,l-1$ do

2. Generate an indicator set $\Omega_{t}$, and choose $\eta m$ rows of $\mathbf{A}$ by using simple random sampling.

3. Compute $\beta_{(jt)}=\frac{|\mathbf{b}_{(jt)}-\mathbf{A}_{(jt)}\mathbf{x}_{t}|}{\|\mathbf{A}_{(jt)}\|_{2}}$ with $jt\in\Omega_{t}$.

4. Sort the $\{\beta_{(jt)}\}$'s and select the $k_{max}$ largest ones. Construct the index set $J_{t}$.

5. Update~$\mathbf{x}_{t+1}=\mathbf{x}_{t}+\mathbf{A}_{J_{t}}^{\dag}(\mathbf{b}_{J_{t}}-\mathbf{A}_{J_{t}}\mathbf{x}_{t})$.

6. If $\mathbf{x}_{t+1}$ is accurate enough, then stop, else continue.

7. End for
\end{algorithm}

The following theorem shows the convergence Algorithm 4, whose proof is similar to that of Algorithm 3. We give the proof for completeness.
\begin{thm}
The iteration sequence $\{\mathbf{x}_{t}\}_{t=0}^{\infty}$  generated by  Algorithm 4 converges to the least-squares solution $\mathbf{x}_{\star}=\mathbf{A}^{\dagger}\mathbf{b}$, with
\begin{equation}
\|\mathbf{x}_{t+1}-\mathbf{x}_{\star}\|_{2}^{2}\leq\left( 1-\hat{\zeta}_{t}\right)\|\mathbf{x}_{t}-\mathbf{x}_{\star}\|_{2}^{2}
\end{equation}
where
$$
\hat{\zeta}_{t}=\frac{\|\mathbf{A}_{J_{t}}\|_{F}^{2}}{\|\mathbf{A}_{t}\|_{F}^{2}+ (\hat{\lambda}_{t}-1)\|\mathbf{A}_{J_{t}}\|_{F}^{2}}\frac{ \sigma_{\min}^{2}(\mathbf{A}_{t}) }{\sigma_{\max}^{2}(\mathbf{A}_{J_{t}})},
$$
and
$$\hat{\lambda}_{t}=\max\limits_{jt\in J_{t}}\left\{\frac{\frac{|\mathbf{b}_{(jt)}-\mathbf{A}_{(jt)}\mathbf{x}_{t}|^{2}}{\|\mathbf{A}_{(jt)}\|_{2}^{2}}}{{\min\limits_{tt\in J_{t}}\left\{\frac{|\mathbf{b}_{(tt)}-\mathbf{A}_{(tt)}\mathbf{x}_{t}|^{2}}{\|\mathbf{A}_{(tt)}\|_{2}^{2}}\right\}}} \right\}\geq 1,\quad jt\in \Omega _{t}.
$$
Here $\mathbf{A}_{t},\mathbf{b}_{t}$ stand for the submatrix and subvector of $\mathbf{A},\mathbf{b}$ corresponding to the
indicator set $\Omega_{t}$.
\end{thm}
\begin{proof}
Similar to \eqref{39}, we have that
\begin{small}\begin{align}\nonumber
   \|\mathbf{A}_{J_{t}}^{\dag}\mathbf{A}_{J_{t}}(\mathbf{x}_{t}-\mathbf{x}_{\star})\|_{2}^{2}
   \geq&\sigma_{\min}^{2}(\mathbf{A}_{J_{t}}^{^{\dag}})\min\limits_{tt\in J_{t}}\left\{\frac{|\mathbf{b}_{(tt)}-\mathbf{A}_{(tt)}\mathbf{x}_{t}|^{2}}{\|\mathbf{A}_{(tt)}\|_{2}^{2}}\right\}\sum\limits_{it\in J_{t}}\|\mathbf{A}_{(it)}\|_{2}^{2}\\ \nonumber
   =&\sigma_{\min}^{2}(\mathbf{A}_{J_{t}}^{^{\dag}})\|\mathbf{A}_{J_{t}}\|_{F}^{2}\min\limits_{tt\in J_{t}}\left\{\frac{|\mathbf{b}_{(tt)}-\mathbf{A}_{(tt)}\mathbf{x}_{t}|^{2}}{\|\mathbf{A}_{(tt)}\|_{2}^{2}}\right\}\\ \nonumber
  \geq &\frac{\|\mathbf{A}_{J_{t}}\|_{F}^{2}}{\sigma_{\max}^{2}(\mathbf{A}_{J_{t}})}\min\limits_{tt\in J_{t}}\left\{\frac{|\mathbf{b}_{(tt)}-\mathbf{A}_{(tt)}\mathbf{x}_{t}|^{2}}{\|\mathbf{A}_{(tt)}\|_{2}^{2}}\right\}\frac{\| \mathbf{A_{t}}\mathbf{x}_{\star}-\mathbf{A_{t}}\mathbf{x}_{t}\|_{2}^{2}}{\| \mathbf{b_{t}}-\mathbf{A_{t}}\mathbf{x}_{t}\|_{2}^{2}}\\ \nonumber
   \geq&\frac{\|\mathbf{A}_{J_{t}}\|_{F}^{2}}{\sigma_{\max}^{2}(\mathbf{A}_{J_{t}})}\sigma_{\min}^{2}(\mathbf{A}_{t})\frac{\min\limits_{tt\in J_{t}}\left\{\frac{|\mathbf{b}_{(tt)}-\mathbf{A}_{(tt)}\mathbf{x}_{t}|^{2}}{\|\mathbf{A}_{(tt)}\|_{2}^{2}}\right\} \|\mathbf{x}_{t}-\mathbf{x}_{\star}\|_{2}^{2}}{\left(\sum\limits_{jt\in J_{t}}\frac{|\mathbf{b}_{(jt)}-\mathbf{A}_{(jt)}\mathbf{x}_{t}|^{2}}{\|\mathbf{A}_{(jt)}\|_{2}^{2}} \|\mathbf{A}_{(jt)}\|_{2}^{2} +\sum\limits_{jt\in [t]/J_{t}}\frac{|\mathbf{b}_{(jt)}-\mathbf{A}_{(jt)}\mathbf{x}_{t}|^{2}}{\|\mathbf{A}_{(jt)}\|_{2}^{2}} \|\mathbf{A}_{(jt)}\|_{2}^{2}\right)} \\ \nonumber
   \geq&\frac{\sigma_{\min}^{2}(\mathbf{A}_{t})}{\sigma_{\max}^{2}(\mathbf{A}_{J_{t}})}\|\mathbf{A}_{J_{t}}\|_{F}^{2}\frac{\|\mathbf{x}_{t}-\mathbf{x}_{\star}\|_{2}^{2}}{\sum\limits_{jt\in J_{t}}\frac{\frac{|\mathbf{b}_{(jt)}-\mathbf{A}_{(jt)}\mathbf{x}_{t}|^{2}}{\|\mathbf{A}_{(jt)}\|_{2}^{2}}}{\min\limits_{tt\in J_{t}}\left\{\frac{|\mathbf{b}_{(tt)}-\mathbf{A}_{(tt)}\mathbf{x}_{t}|^{2}}{\|\mathbf{A}_{(tt)}\|_{2}^{2}}\right\}}\|\mathbf{A}_{(jt)}\|_{2}^{2}+\sum\limits_{jt\in [t]/ J_{t}}\frac{\frac{|\mathbf{b}_{(jt)}-\mathbf{A}_{(jt)}\mathbf{x}_{t}|^{2}}{\|\mathbf{A}_{(jt)}\|_{2}^{2}}}{\min\limits_{tt\in J_{t}}\left\{\frac{|\mathbf{b}_{(tt)}-\mathbf{A}_{(tt)}\mathbf{x}_{t}|^{2}}{\|\mathbf{A}_{(tt)}\|_{2}^{2}}\right\}}\|\mathbf{A}_{(jt)}\|_{2}^{2}},
   \end{align}
\end{small}
where $\mathbf{A}_{[t]/ [J_{t}]}$  is the submatrix obtained from removing $\mathbf{A}_{J_{t}}$
from $\mathbf{A}_{t}$, and the row index set corresponding to $\mathbf{A}_{[t]/ [J_{t}]}$ is denoted by $[t]/ [J_{t}]$.

Moreover, we have that
\begin{align}\nonumber
   \|\mathbf{A}_{J_{t}}^{\dag}\mathbf{A}_{J_{t}}(\mathbf{x}_{t}-\mathbf{x}_{\star})\|_{2}^{2}
   \geq&\frac{\sigma_{\min}^{2}(\mathbf{A}_{t})}{\sigma_{\max}^{2}(\mathbf{A}_{J_{t}})}\|\mathbf{A}_{J_{t}}\|_{F}^{2}\frac{\|\mathbf{x}_{t}-\mathbf{x}_{\star}\|_{2}^{2}}{\sum\limits_{jt\in J_{t}}\hat{\lambda}_{t}\|\mathbf{A}_{(jt)}\|_{2}^{2}+\sum\limits_{jt\in[t]/J_{t}}\|\mathbf{A}_{(jt)}\|_{2}^{2}} \\\nonumber
   \geq&\frac{\sigma_{\min}^{2}(\mathbf{A}_{t})}{\sigma_{\max}^{2}(\mathbf{A}_{J_{t}})}\|\mathbf{A}_{J_{t}}\|_{F}^{2}\frac{\|\mathbf{x}_{t}-\mathbf{x}_{\star}\|_{2}^{2}}{ (\hat{\lambda}_{t}-1)\|\mathbf{A}_{J_{t}}\|_{F}^{2}+\|\mathbf{A}_{t}\|_{F}^{2}}\\
   =&\frac{\|\mathbf{A}_{J_{t}}\|_{F}^{2}}{\|\mathbf{A}_{t}\|_{F}^{2}+ (\hat{\lambda}_{t}-1)\|\mathbf{A}_{J_{t}}\|_{F}^{2}}\frac{ \sigma_{\min}^{2}(\mathbf{A}_{t}) }{\sigma_{\max}^{2}(\mathbf{A}_{J_{t}})}\|\mathbf{x}_{t}-\mathbf{x}_{\star}\|_{2}^{2}.\label{313}
\end{align}

%

Substituting \eqref{313} into \eqref{36}, we get
\begin{equation*}
~~~~~~\|\mathbf{x}_{t+1}-\mathbf{x}_{\star}\|_{2}^{2}\leq\left( 1-\frac{\|\mathbf{A}_{J_{t}}\|_{F}^{2}}{\|\mathbf{A}_{t}\|_{F}^{2}+ (\hat{\lambda}_{t}-1)\|\mathbf{A}_{J_{t}}\|_{F}^{2}}\frac{ \sigma_{\min}^{2}(\mathbf{A}_{t}) }{\sigma_{\max}^{2}(\mathbf{A}_{J_{t}})}\right)\|\mathbf{x}_{t}-\mathbf{x}_{\star}\|_{2}^{2},
\end{equation*}
which completes the proof.
\end{proof}

\subsection{Solving large-scale linear system with {\it multiple} right-hand sides}
In this subsection, we focus on the large-scale linear system \eqref{1.1} with multiple right-hand sides. As far as we are aware,  
there are few Kaczmarz-type methods for solving this problem. We will propose a semi-randomized {\it block} Kaczmarz method with simple random sampling for
large-scale linear system with multiple right-hand sides.


Inspired by Algorithm 3, for the $j$-th linear system, 
we select the working row $\tilde{t}_{j}$ according to
\begin{equation}
\max\limits_{1 \leq t_{i}\leq m}\left\{\frac{|\mathbf{B}_{(t_{i},j)}-\mathbf{A}_{(t_{i})}\mathbf{x}_{j}^{(t)}|}{\|\mathbf{ A}_{(t_{i})} \|_{2}}\right\},\quad j=1,2,\ldots,k_{b}, 
\end{equation}
and update the $j$-th solution vector ${\bf x}_j$ in the following way
\begin{equation}\label{315}
\mathbf{x}_{j}^{(t+1)}=\mathbf{x}_{j}^{(t)}+\frac{\mathbf{B}_{(\widetilde{t}_{j},j)}-\mathbf{A}_{(\widetilde{t}_{j})}\mathbf{x}_{j}^{(t)} }{\|\mathbf{A}_{(\widetilde{t}_{j})}\|_{2}^{2}}\mathbf{A}_{(\widetilde{t}_{j})}^{T},\quad j=1,2,\ldots,k_{b}.
\end{equation}

Let
$$\mathbf{X}^{(t+1)}=
\begin{pmatrix}
\mathbf{x}_{1}^{(t+1)}, & \mathbf{x}_{2}^{(t+1)}, & \ldots, & \mathbf{x}_{k_{b}}^{(t+1)}
\end{pmatrix},\quad t=0,1,\ldots
$$
be the approximation of ${\bf X}$ in the $(t+1)$-th iteration,
then \eqref{315} can be written as
$$
\mathbf{X}^{(t+1)}=\mathbf{X}^{(t)}+
\begin{pmatrix}
\mathbf{A}^{T}_{(\widetilde{t}_{1})},
\ldots,
\mathbf{A}^{T}_{(\widetilde{t}_{k_{b}})}
\end{pmatrix}^{T}
diag\left(\frac{\mathbf{B}_{(\widetilde{t}_{1},1)}-\mathbf{A}_{(\widetilde{t}_{1})}\mathbf{x}_{1}^{(t)}}{\| \mathbf{A}_{\widetilde{t}_{1}}\|_{2}^{2}},\ldots,\frac{\mathbf{B}_{(\widetilde{t}_{k_{b}},k_{b})}-\mathbf{A}_{(\widetilde{t}_{k_{b}})}x_{k_{b}}^{(t)}}{\| \mathbf{A}_{(\widetilde{t}_{k_{b}})}\|_{2}^{2}}\right).
$$

We are in a position to present the following randomized block Kaczmarz method with simple random sampling. Notice that one can update all the 
solutions simultaneously.
\begin{algorithm}
\caption{A semi-randomized {\it block} Kaczmarz method with simple random sampling for
large-scale linear system with {\it multiple} right-hand sides}


${\rm \bf Input}$:

\ \  $\mathbf{A},\ \mathbf{B},\ l$,  $\mathbf{X}_{0}$, $0<\eta\ll 1$.

{\rm \bf Output}:

\ \  $\mathbf{X}^{(l)}$ .

1. For $t=0,1,\ldots,l-1$ do

2. Generate an indicator set $\Omega_{t}$, and choose $\eta m$ rows of $A$ by using the simple random sampling method.

3. Select the working row number~$\widetilde{t}_{j}\in \Omega_{t}$ of the first~$j$-th righthand side according to~$\max\limits_{t_{i}\in \Omega_{t}}\left\{\frac{|\mathbf{B}_{(t_{i},j)}-\mathbf{A}_{(t_{i})}\mathbf{x}_{j}^{(t)}|}{\|\mathbf{ A}_{(t_{i})} \|_{2}}\right\}$,~$j=1,2,\ldots,k_{b}$.

4. Update\\
$\mathbf{X}^{(t+1)}=\mathbf{X}^{(t)}+
\begin{pmatrix}
\mathbf{A}^{T}_{(\widetilde{t}_{1})},
\ldots,
\mathbf{A}^{T}_{(\widetilde{t}_{k_{b}})}
\end{pmatrix}^{T}
diag\left(\frac{\mathbf{B}_{(\widetilde{t}_{1},1)}-\mathbf{A}_{(\widetilde{t}_{1})}\mathbf{x}_{1}^{(t)}}{\| \mathbf{A}_{\widetilde{t}_{1}}\|_{2}^{2}},\ldots,\frac{\mathbf{B}_{(\widetilde{t}_{k_{b}},k_{b})}-\mathbf{A}_{(\widetilde{t}_{k_{b}})}x_{k_{b}}^{(t)}}{\| \mathbf{A}_{\widetilde{t}_{k_{b}}}\|_{2}^{2}}\right)
$
5. If~$\mathbf{X}^{t+1}$ is accurate enough, then stop, else continue.

6. End for
\end{algorithm}

The following theorem shows the convergence of Algorithm 5.

\begin{thm}
The iteration sequence $\{\mathbf{X}^{(t)}\}_{t=0}^{\infty}$  generated by Algorithm 5 converges to the least-squares solution $\mathbf{X}^{\star}=\mathbf{A}^{\dagger}\mathbf{B}$, with
\begin{equation}
\mathbb{E}\|\mathbf{X}^{(t+1)}-\mathbf{X}^{\star}\|_{F}^{2}\leq\left( 1-\frac{\sigma_{\min}^{2}(\mathbf{A}_{t})}{\|\mathbf{A}_{t}\|_{F}^{2}}\right)\|\mathbf{X}^{(t)}-\mathbf{X}^{\star}\|_{F}^{2},
\end{equation}
where $\mathbf{A}_{t},\mathbf{B}_{t}$ represent the submatrices of $\mathbf{A}$, $\mathbf{B}$ corresponding to the
indicator set $\Omega_{t}$, respectively.
\end{thm}
\begin{proof}
Let $\mathbf{X}^{\star}=\mathbf{A}^{\dagger}\mathbf{B}$
be the least-squares solution of \eqref{1.1}, then
$\mathbf{x}_{j}^{\star}=\mathbf{A}^{\dagger}\mathbf{B}_{(:,j)}$, $j=1,2,\ldots,k_{b}$,
where $\mathbf{B}_{(:,j)}$ is the $j$-th column of $\mathbf{B}$.
Furthermore,
\begin{equation}
\mathbb{E}\|\mathbf{X}^{(t+1)}-\mathbf{X}^{\star}\|_{F}^{2}=\mathbb{E}\sum\limits_{j=1}\limits^{k_{b}}\|\mathbf{x}_{j}^{(t+1)}-\mathbf{x}_{j}^{\star}\|_{2}^{2}
\\=\sum\limits_{j=1}\limits^{k_{b}}\mathbb{E}\|\mathbf{x}_{j}^{(t+1)}-\mathbf{x}_{j}^{\star}\|_{2}^{2}.\label{317}
\end{equation}
We notice that
$$
\mathbb{E}\|\mathbf{x}_{j}^{(t+1)}-\mathbf{x}_{j}^{\star}\|_{2}^{2}=\|\mathbf{x}_{j}^{(t)}-\mathbf{x}_{j}^{\star}\|_{2}
^{2}-\mathbb{E}\|\mathbf{x}_{j}^{(t+1)}-\mathbf{x}_{j}^{(t)}\|_{2}^{2}, \quad j=1,2,\ldots,k_{b},
$$
and
$$\mathbb{E}\|\mathbf{x}_{j}^{(t+1)}-\mathbf{x}_{j}^{(t)}\|_{2}^{2}
    =\frac{|\mathbf{B}_{(\widetilde{t}_{j},j)}-\mathbf{A}_{(\widetilde{t}_{j})}\mathbf{x}_{j}^{(t)} |^{2}}{\|\mathbf{A}_{(\widetilde{t}_{j})}\|_{2}^{2}}
   =\max\limits_{1\leq i\leq \eta m}\left\{\frac{|\mathbf{B}_{(t_{i},j)}-\mathbf{A}_{(t_{i})}\mathbf{x}_{j}^{(t)} |^{2}}{\|\mathbf{A}_{(t_{i})}\|_{2}^{2}}\right\}.
$$

Since
\begin{small}\begin{align*}\nonumber
   \max\limits_{1\leq i\leq \eta m}\left\{\frac{|\mathbf{B}_{(t_{i},j)}-\mathbf{A}_{(t_{i})}\mathbf{x}_{j}^{(t)} |^{2}}{\|\mathbf{A}_{(t_{i})}\|_{2}^{2}}\right\}
=&\frac{\max\limits_{1\leq i\leq \eta m}\left\{\frac{|\mathbf{B}_{(t_{i},j)}-\mathbf{A}_{(t_{i})}\mathbf{x}_{j}^{(t)} |^{2}}{\|\mathbf{\mathbf{A}}_{(t_{i})}\|_{2}^{2}}\right\}}{\|\mathbf{\mathbf{A}}_{t}\mathbf{x}_{j}^{(t)}-\mathbf{B}_{(:,j)}\|_{2}^{2}}\|\mathbf{A}_{t}\mathbf{x}_{j}^{(t)}-\mathbf{B}_{(:,j)}\|_{2}^{2}\\ \nonumber
   =&\frac{\max\limits_{1\leq i\leq \eta m}\left\{\frac{|\mathbf{B}_{(t_{i},j)}-\mathbf{A}_{(t_{i})}\mathbf{x}_{j}^{(t)} |^{2}}{\|\mathbf{A}_{(t_{i})}\|_{2}^{2}}\right\}\|\mathbf{A}_{t}\mathbf{x}_{j}^{(t)}-\mathbf{B}_{(:,j)}\|_{2}^{2}}{\sum\limits_{i=1}\limits^{\eta m}\|\mathbf{A}_{(t_{i})}\|_{2}^{2}\frac{|\mathbf{B}_{((t_{i},j)}-\mathbf{A}_{(t_{i})}\mathbf{x}_{j}^{(t)}|^{2}}{\|\mathbf{A}_{(t_{i})}\|_{2}^{2}}}\\ \nonumber
   \geq&\frac{\|\mathbf{A}_{t}\mathbf{x}_{j}^{(t)}-\mathbf{B}_{(:,j)}\|_{2}^{2}}{\sum\limits_{i=1}\limits^{\eta m}\|\mathbf{A}_{(t_{i})}\|_{2}^{2}}\\
   \geq&\frac{\sigma_{\min}^{2}(\mathbf{A}_{t})}{\|\mathbf{A}_{t}\|_{F}^{2}}\|\mathbf{x}_{j}^{(t)}-\mathbf{x}_{j}^{\star}\|_{2}^{2},
\end{align*}
\end{small}
we have that
\begin{align}\nonumber
\mathbb{E}\|\mathbf{x}_{j}^{(t+1)}-\mathbf{x}_{j}^{\star}\|_{2}^{2}
=&\|\mathbf{x}_{j}^{(t)}-\mathbf{x}_{j}^{\star}\|_{2}^{2}-\mathbb{E}\|\mathbf{x}_{j}^{(t+1)}-\mathbf{x}_{j}^{(t)}\|_{2}^{2} \\ \leq&\left (1-\frac{\sigma_{\min}^{2}(\mathbf{A}_{t})}{\|\mathbf{A}_{t}\|_{F}^{2}}\right)\|\mathbf{x}_{j}^{(t)}-\mathbf{x}_{j}^{\star}\|_{2}^{2}.\label{318}
\end{align}
A combination of \eqref{317} and \eqref{318} gives the result.
\end{proof}




\section{Numerical Results}

In this section, we perform some numerical experiments to show the numerical behavior of the proposed algorithms. All the numerical experiments are run on a Hp workstation with 20 cores double Intel(R)Xeon(R) E5-2640 v3 processors, with CPU 2.60 GHz and RAM 256 GB. 
The operation system is 64-bit Windows 10. The numerical results are obtained from using MATLAB 2018b. In the tables below, we denote by {\tt IT}
the number of iterations, by {\tt CPU} the computing time in seconds. If the number of iteration exceeds $10^6$, or the CPU time exceeds 12 hours, we will stop the algorithm and declare it fails to converge. Thus, in the tables, ``$--$" implies that the number of iterations exceeds $10^6$ or the CPU time exceeds 12 hours, and ``{\tt O.M.}" stands for the algorithm suffers from out of memory. All the experiments are repeated for 5 times, and the iteration numbers as well as the CPU time in seconds are the mean from the 5 runs.

\subsection{Numerical experiments for large-scale linear systems with {\it single} right-hand side}
In this subsection, we are interested in the large-scale linear system \eqref{1.2} with single right-hand side.
In order to show the efficiency of our proposed algorithms Algorithm 3 and Algorithm 4 , we compare them with the following randomized block Kaczmarz methods for large-scale linear systems with single righthand side or multiple righthand sides:\\
~\\
$\mathbf{GBK}$: The greedy block Kaczmarz method due to Niu and Zheng \cite{21}.\\
$\mathbf{GRBK}$: The greedy randomized block Kaczmarz  method due to Liu and Gu \cite{16} .\\
$\mathbf{BEM}$: The block extension of Motzkin method due to Liu and Gu \cite{16} .\\
$\mathbf{RBK}$: The randomized block Kaczmarz method due to Needell and Tropp \cite{20}.\\
$\mathbf{RABK}$: The randomized average block Kaczmarz methoddue to Necoara \cite{18}.\\
$\mathbf{REABK}$: The randomized extended average block Kaczmarz method due to Du, Shi and Sun \cite{5}.\\

In practical implementations, we make use of the same random partition $\Gamma=\{\tau_{1},\tau_{2},\ldots,\tau_{p}\}$ for RBK, GRBK and REABK methods. More precisely, let $N_{r}$ be the number of rows for each block and $p=\lceil \frac{m}{N_{r}}\rceil$. The row blocks of the RBK , GRBK and REABK algorithm are uniformly divided into $p$ blocks:
$$
\tau_{i}=\left\{(i-1)N_{r}+1,(i-1)N_{r}+2,\ldots,iN_{r}\right\},~~i=1,2,\ldots,p-1,
$$
and
$$
\tau_{p}=\left\{(p-1)N_{r}+1,(p-1)N_{r}+2,\ldots,m\right\},~N_{\tau_{p}}\leq N_{r}.
$$

In the GBK method, the parameter $\varepsilon$ used in the greedy strategy is chosen as \cite{21}
$$
\varepsilon=\frac{1}{2}+\frac{1}{2}\frac{\|\mathbf{b}-\mathbf{Ax}_{t}\|_{2}^{2}}{\|\mathbf{A}\|_{F}^{2}}\left(\max\limits_{1\leq i\leq m}\left(\frac{|\mathbf{b}_{(i)}-\mathbf{A}_{(i)}\mathbf{x}_{t}|^{2}}{\|\mathbf{A}_{(i)}\|_{2}^{2}}\right)\right)^{-1}.
$$
In the RBK, GRBK and REABK methods, a fixed row block size $N_{r}$ and a fixed column block size $N_{c}$ are needed. For the sake of fairness, the row block size $N_{r}$ and column block size  $N_{c}$ are set to be the same as the value of $k_{max}$ initialized in Algorithm 4. That is,
$$
N_r=N_c=k_{max}.
$$

In the REABK method, the empirical step size $\alpha$ is chosen in the interval $(\frac{1}{\beta_{max}},\frac{2\|\mathbf{A}\|_{F}^{2}}{\sigma^{2}_{max}(\mathbf{A})}) $ \cite{5}. In the BEM method, as the relaxation parameter $\rho$ is taken from the set $\{0.1,0.2,\ldots,0.9\}$ \cite{16}, we choose it to be 0.5 in all the experiments. In the RABK algorithm, as the step size $\alpha_{k}\in (0, 2)$ \cite{18}, we choose it as 1.95.

The right-hand side ~${\bf b} \in R^{m}$ is taken to be $\mathbf{Ax}^{\star}$, where $\mathbf{x}^{\star}\in R^{n}$ is generated by using the MATLAB built-in function $randn(n,1)$ or $ones(n,1)$. The initial vector ~$\mathbf{x}_{0}$ is taken as the zero vector. The
stopping criterion is
$$
RES=\frac{\|\mathbf{x}_{t}-\mathbf{x}_{\star}\|_{2}^{2}}{\|\mathbf{x}_{\star}\|_{2}^{2}} < tol,
$$
where $\mathbf{x}_{t}$ is a computed solution in the $t$-th iteration, and $tol$ is a user-described tolerance.

\begin{table}[H]
 \begin{small}
 \begin{center}
 \caption{$\mathbf{Example ~4.1}$. Test matrices used in Example 4.1, where ``$(\cdot)^T$'' stands for the transpose of a matrix.}
   \vspace{0.2cm}
\begin{tabular}{|c|c|c|c|c}
\hline  Name&size $(m\times n)$&nnz&background\\
\hline lp\_pds\_$10^{T}$&$49932 \times 16558$&853829&Linear Programming Problem\\
\hline mesh\_deform&$234032\times 9393$&853829&Coumputer Graphics/Vision Problem\\
\hline sls&$1748122\times 62729$&6804304&Least Squares Problem\\
\hline cage10&$11397\times11397$&150645&Background Directed Weighted Graph\\
\hline model10$^{T}$&$16819\times4400$&150372&Linear Programming Problem\\
\hline Alemdar&$6245\times6245$&42581&2D/3D Problem\\
\hline rajat13&$7598\times7598$&48762&Circuit Simulation Problem\\
\hline nemsemm2$^{T}$&$48873\times6943$&182010&Linear Programming Problem\\
\hline rlfddd$^{T}$&$61521\times4050$&264627&Linear Programming Problem\\
\hline lp\_pds\_06$^{T}$&$29351\times9881$&63220&Linear Programming Problem\\
\hline
\end{tabular}
\end{center}
\end{small}
\end{table}

$\mathbf{Example ~4.1}$. 
In this example, the test matrices are from the University of Florida sparse matrix collection \cite{6}. Table 4.1 presents the details of these matrices. We run Algorithm 3, Algorithm 4 and the six randomized block Kaczmarz methods for the ten test problems. In this example, $\mathbf{x}^{\star}\in R^{n}$ is generated by using $randn(n,1)$.  Numerical results are shown in Table 4.2.

Some remarks are in order. First, it is observed in Table 4.2 that both Algorithm 3 and Algorithm 4 run much faster than the other randomized block Kaczmarz methods. Based on the technique of simple random sampling and the idea of semi-randomization, Algorithm 4 makes use of a small section of rows of $A$ for choosing the working block. Consequently, the computational cost in each iteration is much less than those of the others, and it often performs the best in terms of CPU time.
Second, some methods such as GRBK, RBK, RABK and REABK fail to converge within the specified time or number of iterations, while the two proposed methods work quite well. Specifically, we find that REABK does not work for all the problems. This is because one has to determine many parameters in advance, which are difficult to tune.

Third, we see that Algorithm 3 and Algorithm 4 often use much more iterations than the other methods such as GBK and GBM. This is because the cost used in each iteration of our methods is much less than that of the others. 
Two show this more precisely, for the {\tt mesh\_deform} matrix, we plot in Figure 1 the number of rows in the row index set of the GBK (left) and the BEM (right) methods during iterations. As the working sets of GBK and BEM are constructed during iterations, it is seen from the figure that
the block sizes of the two methods can be very large. 
More precisely, one sees that the maximum value of $Nu(\mathbf{A}_{J_{t}})$ is above 4000 for the GBK method, and is 3500 above for the BEM method. As a comparison, in Algorithm 3 and Algorithm 4, we have $(Nu(\mathbf{A}_{J_{t}})=k_{max})$, which is fixed. Thus, Algorithm 4 effectively avoids the disadvantage of GBK and BEM algorithm. That is,  $Nu(\mathbf{A}_{J_{t}})$  can be too large in these two methods,  which will lead to great computational overhead for large-scale matrices.



\begin{table}[H]
\begin{small}
  \begin{center}
  \caption
  {$\mathbf{Example ~4.1}$. Numerical results of the randomized block Kaczmarz-type methods for solving \eqref{1.2} on some large-sparse matrices, $tol=10^{-3}$, where the best ones in terms of CPU time are in bold. The sampling ratio is chosen as $\eta=0.1$ in Algorithm 4, and $k_{max}=10$.}
  \vspace{0.2cm}
    \begin{tabular}{lrrrrrrrrr}
     \hline
      Name&         &Alogrithm 3  &Alogrithm 4&GBK   & GRBK  & BEM     & RBK & RABK&REABK \\ \hline
  lp\_pds\_$10^{T}$& CPU&228.57 &$\mathbf{152.04}$    &7019.40 &$--$&5870.20&1496.8&1186.70 &$--$        \\
             & IT     &26873&13964     &681  &$--$ &679     &205392 &972419 &$--$ \\ \hline
   mesh\_deform  & CPU &90.86&$\mathbf{78.96}$ &223.73    &$--$&538.85 &419.36 &2827.30  &$--$             \\
                & IT  &3698&2885     &84     &$--$&700     &572.63  &766772 &$--$\\ \hline
sls            & CPU  &$\mathbf{2265.90}$&5796.20   &32309.00    &$O.M.$  &$O.M.$  &22155.00  &$--$&$--$ \\
               & IT   &1606&21450     &12   &$O.M.$      &$O.M.$        &604617  &$--$ &$--$\\ \hline
     cage10    &  CPU&$\mathbf{7.12}$ & 9.07    &11.63  &12.97   &47.36 &47.18   &16.72  &$--$              \\
              & IT   & 1606&1613      &48      &2311   &36    &15377   &61621  &$--$   \\ \hline
        model$10^{T}$ &CPU& 126.87&83.34    &$\mathbf{66.33}$  &$--$   &507.13 &$--$   &$--$ &$--$             \\
              & IT   &41582&21024     &1742 &$--$     &32642   &$--$ &$--$  &$--$  \\ \hline
        Alemdar   &  CPU &$\mathbf{546.32}$ &  686.07  &6472.30  &2052.80  &4697.40 &2309.10   &$--$  &$--$            \\
              & IT &211792  &171151    &20882      &39406    &18820    &62992    &$--$ &$--$  \\ \hline
       rajat13   &  CPU &32.66&$\mathbf{31.55 }$     &58.64 &1576.70 &$--$ &$--$   &$--$  &$--$             \\
             & IT & 14710 &11729  &2073     &14706   &$--$   &$--$  &$--$  &$--$   \\ \hline
        nemsemm2$^{T}$  &CPU&$\mathbf{3.87 }$ &6.45      &7.97 &321.68  &69.71 &$--$  &$--$  &$--$             \\
              & IT  &874 &809&36 &990   &2647   &$--$   &$--$  &$--$  \\ \hline
        rlfddd$^{T}$ &CPU & 3.41&$\mathbf{2.53}$     &5.16 &3529.80  &11.51 &37.28  &31.70 &$--$            \\
              & IT   &660&396&16 &1104   &169   &14893  &48231 &$--$  \\ \hline
        lp\_pds\_06$^{T}$& CPU& 37.01&$\mathbf{32.12}$     &324.24 &10776.00  &321.13 &369.17 &306.19 &$--$             \\
              & IT  &8494 &6082&449&15124  &392   &104036&482416 &$--$  \\ \hline

\end{tabular}
\end{center}
\end{small}
\end{table}

\begin{figure}[H]
\begin{center}
{\centering\includegraphics[width=5cm]{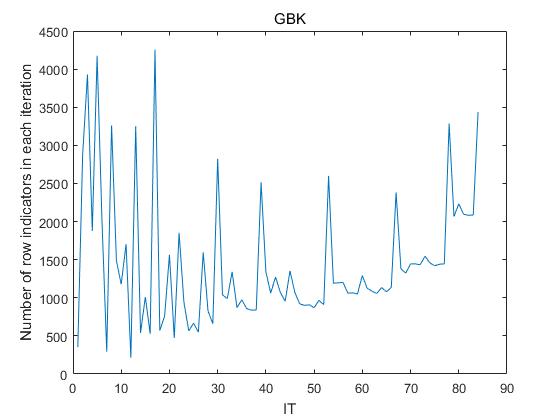}}{\centering\includegraphics[width=5cm]{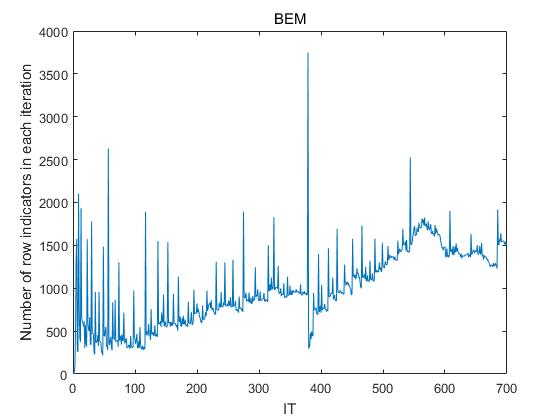}}\\

{\small {\bf Figure 1.} The number of rows $Nu(\mathbf{A}_{J_{t}})$ in the row index set of the GBK (left) and the BEM (right) methods in each iteration. The  matrix {\tt mesh\_deform} matrix.
}
\end{center}
\end{figure}

$\mathbf{Example~4.2}$. In this example, the matrices are two large-dense matrices {\tt RCVI\_4class} (of size 29992$\times$9625) and {\tt 20NewSHomeSort} (of size 61188$\times$18744), which are text data \footnote{http://www.cad.zju.edu.cn/home/dengcai}. In this example, the ``exact" solution  $\mathbf{x}^{\star}\in R^{n}$ is set to be $ones(n,1)$. We run the eight algorithms on the two large-dense matrices, and the convergence tolerance  $tol$  is chosen as $10^{-3}$. The sampling ratio are chosen as $\eta=0.1$ in Algorithm 4, and we use $k_{max}=10,50,100$, respectively. Table 4.3 lists the numerical results.

\begin{table}
\begin{scriptsize}
\begin{center}
\caption
  {$\mathbf{Example~4.2}$. Numerical results of the randomized block Kaczmarz-type methods for solving \eqref{1.2} on some large-dense matrices, $tol=10^{-3}$, where the best ones in terms of CPU time are in bold. $tol=10^{-3}$. The sampling ratio is chosen as $\eta=0.1$ in Algorithm 4, and $k_{max}=10,50,100$, respectively.}
\begin{tabular}{c|c|c|c|c|c|c|c|c|c|c}
	\hline
    Name&& &Alogrithm 3   &Alogrithm 4&GBK   & GRBK  & BEM     & RBK & RABK&REABK \\ \hline
    \multirow{6}*{RCVI\_4class}
    &\multirow{2}*{$k_{max}=10$}&CPU & $\mathbf{23.74}$ &37.35 &42.20&$--$&65.97&182.54&479.76&$--$\\
    ~&~&IT &2407   &1959 &239 &$--$&1186&38609&97917&$--$\\ \cline{2-11}
    ~&\multirow{2}*{$k_{max}=50$}&CPU& $\mathbf{12.85}$ & 15.21&42.42 &2575.70&66.05&181.86&150.73&$--$\\
    ~&~&IT  &  579   &608 &239 &2600&1186&8979&69627&$--$\\ \cline{2-11}
    ~&\multirow{2}*{$k_{max}=100$}&CPU& 16.66&$\mathbf{13.89}$  &42.33&816.43&66.17&315.32&164.98&$--$\\
    ~&~&IT  &  390   &269  &239 &1608&1186&7640&70112&$--$\\ \hline
   \multirow{6}*{20NewSHomeSort}
    &\multirow{2}*{$k_{max}=10$}&CPU &$\mathbf{4065.10}$ &7099.90 & 16563.00&$--$&$--$&$--$&$--$&$--$\\
    ~&~&IT  & 169976&134026 &18730&$--$&$--$ &$--$&$--$&$--$\\ \cline{2-11}
    ~&\multirow{2}*{$k_{max}=50$}&CPU  & $\mathbf{2489.00}$&2594.30 & 15824.00&$--$&$--$&$--$&$--$&$--$\\
    ~&~&IT   &  67290  &134556&18730&$--$&$--$ &$--$&$--$&$--$\\ \cline{2-11}
    ~&\multirow{2}*{$k_{max}=100$}&CPU  &2555.00 & $\mathbf{1797.80}$& 15916.00&$--$&$--$ &$--$&$--$&$--$\\
    ~&~&IT &   36379   &19658&18730&$--$&$--$ &$--$&$--$&$--$\\ \hline

\end{tabular}
\end{center}
\end{scriptsize}
\end{table}
It is seen from Table 4.3 that the GRBK, BEM, RBK, RABK and REABK methods do not work for the
20NewSHomeSort matrix, which is of size 61188$\times$18744. Moreover, both Algorithm 3 and Algorithm 4 outperforms GBK in terms of CPU time.
For the  RCVI\_4class matrix, Algorithm 3 and Algorithm 4 are about 2--3 times faster than GBK, and are about 10 times faster than RABK. Interestingly, we ee in the Table 4.3 that a larger $k_{max}$ may improve the numerical performances of the proposed algorithms. Thus, the proposed randomized block Kaczmarz methods are more suitable to large-scale and dense linear systems.

\subsection{Numerical experiments for large-scale linear systems with {\it multiple} right-hand sides}

In this subsection, we consider the large-scale linear system \eqref{1.1}: $\mathbf{AX}=\mathbf{B}$, where $\mathbf{B}\in R^{m\times k_{b}}$ is of full column rank. So as to illustrate the efficiency of Algorithm 5, we compare it with three randomized Kaczmarz methods and three randomized block Kaczmarz methods:\\ 
~\\
$\mathbf{GBK}$: The greedy block Kaczmarz method due to Niu and Zheng \cite{21}.\\
$\mathbf{SRK}$: The semi-randomized Kaczmarz method with simple random sampling for
large-scale linear system with {\it simple} right-hand side \cite{12}.\\
$\mathbf{RK}$: The randomized Kaczmarz method proposed in \cite{24}.\\
$\mathbf{GRK}$: The greedy randomized Kaczmarz method due to Bai and Wu \cite{3}.\\
$\mathbf{REABK}$: The randomized extended average block Kaczmarz method  proposed by Du, Shi and Sun \cite{5}.\\
$\mathbf{RBK}$: The randomized block Kaczmarz method due to Needell and Tropp \cite{20}.\\

In the GBK, RBK and REABK methods, the row block size $N_{r}$ is set to be 100.
As the six compared methods are just for single right-hand side linear systems, we have to solve the $k_{b}$ linear systems one-by-one using these methods. 
As a comparison, we can update all the solutions of large-scale linear system with multiple right-hand sides simultaneously in Algorithm 5.
Thus, except for Algorithm 5, the value of ``IT" is rounded up by the average number of solving the $k_{b}$ linear systems with single right-hand side.

In this subsection, the right-hand sides $\mathbf{B}\in R^{m\times k_{b}}$ is taken to be $\mathbf{AX}^{\star}$, where $\mathbf{X}^{\star}\in R^{n\times k_{b}}$ is generated by using the MATLAB built-in function $randn(n,k_{b})$. The initial matrix $\mathbf{X}_{0}$ is chosen as the zero matrix. Let $\mathbf{X}^{\star}=
\begin{pmatrix}
\mathbf{x}_{1}^{\star}, & \mathbf{x}_{2}^{\star}, & \ldots, & \mathbf{x}_{k_{b}}^{\star}
\end{pmatrix}$
and
$\mathbf{X}^{(t)}=
\begin{pmatrix}
\mathbf{x}_{1}^{(t)}, & \mathbf{x}_{2}^{(t)}, & \ldots, & \mathbf{x}_{k_{b}}^{(t)}
\end{pmatrix}$ 
be a computed solution in the $t$-th iteration, then the
stopping criterion is
$$
RES=\max\limits_{1\leq j\leq k_{b}}\left\{\frac{\|\mathbf{x}_{j}^{(t)}-\mathbf{x}_{j}^{\star}\|_{2}^{2}}{\|\mathbf{x}_{j}^{\star}\|_{2}^{2}}\right\} < tol,
$$
where $tol$ is a user-described tolerance. 
\begin{table}
\begin{scriptsize}
\begin{center}
\caption
  {$\mathbf{Example~4.3}$. Numerical results of the methods on linear systems whose coefficient matrix $\mathbf{A}\in R^{m\times n}$ is generated by using the MATLAB build-in function randn(m,n), with $n=500$, $m=5000,6000,7000$, tol=$10^{-6}$. The sampling ratio is chosen as $\eta=0.01$ in Algorithm 5. Except for Algorithm 5, the value of ``IT" is rounded up by the average number of solving the $k_{b}$ linear systems with single right-hand side.}
\begin{tabular}{c|c|c|c|c|c|c|c|c|c}
	\hline
    Name&&&${\rm Algorithm 5}$&GBK&SRK&RK&GRK&REABK&RBK\\ \hline
      \multirow{6}*{randn(5000,500)}
    &\multirow{2}*{$k_{b}=10$}&CPU &$\mathbf{0.72}$  & 1.22&3.41&16.54&7.04&$--$&4.02\\
    ~&~&IT    &1251 &32 &1239&7759&955&$--$&71\\ \cline{2-10}
    ~&\multirow{2}*{$k_{b}=20$}&CPU&$\mathbf{1.13}$  &2.40 &6.54&33.16&13.62&$--$&8.07\\
    ~&~&IT       &1259 &32 &1231&7740&943&$--$&71\\ \cline{2-10}
    ~&\multirow{2}*{$k_{b}=50$}&CPU&$\mathbf{1.72}$  &5.88 &16.96&82.33&34.08&$--$&20.11\\
    ~&~&IT       &1271  &32 &1231&7774&952&$--$&70\\ \hline
   \multirow{6}*{randn(6000,500)}
    &\multirow{2}*{$k_{b}=10$}&CPU & $\mathbf{0.81}$& 1.24&3.54&16.61 &7.33&$--$&3.84\\
    ~&~&IT    &1160&28 &1145&7596&898&$--$&69\\ \cline{2-10}
    ~&\multirow{2}*{$k_{b}=20$}&CPU&$\mathbf{1.14}$  &2.35 &6.83&33.87&14.19&$--$&7.89\\
    ~&~&IT       &1170 &29 &1148&7562&889&$--$&70\\\cline{2-10}
    ~&\multirow{2}*{$k_{b}=50$}&CPU&$\mathbf{1.65}$  &5.77 &17.80&84.86&35.63&$--$&19.26\\
    ~&~&IT       &1206  &29 &1146&7569&889&$--$&69\\ \hline
    \multirow{6}*{randn(7000,500)}
    &\multirow{2}*{$k_{b}=10$}&CPU &$\mathbf{0.86}$  &1.18&4.14&18.85&8.02&$--$&3.90\\
    ~&~&IT    &1123 &26&1090&7483&849&$--$&69\\ \cline{2-10}
    ~&\multirow{2}*{$k_{b}=20$}&CPU&$\mathbf{1.15}$  &2.32&7.98&37.17&15.71&$--$&7.61\\
    ~&~&IT       &1136 &26 &1089&7451&849&$--$&67\\ \cline{2-10}
    ~&\multirow{2}*{$k_{b}=50$}&CPU&$\mathbf{1.67}$  &5.99 &21.67&93.67&38.92&$--$&19.19\\
    ~&~&IT       &1125 &27 &1089&7451&849&$--$&68\\ \hline

\end{tabular}
\end{center}
\end{scriptsize}
\end{table}

$\mathbf{Example~4.3}$. In this example, the test matrix $A$ is generated by using $randn(m,n)$, with $n=500$, $m=5000,6000,7000$ for the first test, and $m=4000$, $n=500,600,700$ for the second test. The right-hand side $\mathbf{B}=\mathbf{AX}^{\star}\in R^{m\times k_{b}}$, where $\mathbf{X}^{\star}\in R^{n\times k_{b}}$ is generated by using the MATLAB function $randn(n,k_{b})$. Numerical results are given in Table 4.4 and Table 4.5. We see from the numerical results of the two tables that Algorithm 5 performs the best in terms of CPU time, especially when $k_{b}$ is relatively large. Although Algorithm 5 may use much more iterations than the others, it converges the fastest. The reason is that there is no need to compute all the elements of the residual vectors, nor scanning the whole data matrix. Notice that REABK still do not work in this example.  

\begin{table}
\begin{scriptsize}
\begin{center}
\caption
  {$\mathbf{Example~4.3}$. Numerical results of the methods on linear systems whose coefficient matrix $\mathbf{A}\in R^{m\times n}$ is generated by using the MATLAB build-in function randn(m,n), $m=4000$ and $n=500,600,700$, $tol=10^{-6}$. The sampling ratio is chosen as $\eta=0.01$ in Algorithm 5. Except for Algorithm 5, the value of ``IT" is rounded up by the average number of solving the $k_{b}$ linear systems with single right-hand side.}
\begin{tabular}{c|c|c|c|c|c|c|c|c|c}
	\hline
    Name&&&${\rm Algorithm~5}$&GBK&SRK&RK&GRK&REABK&RBK\\ \hline
       \multirow{6}*{randn(4000,500)}
    &\multirow{2}*{$k_{b}=10$}&CPU &$\mathbf{0.66}$ &2.68&14.39&6.04&3.95&$--$&4.06\\
    &~&IT    &1396&36 &1368&8007&1040&$--$&72\\ \cline{2-10}
    ~&\multirow{2}*{$k_{b}=20$}&CPU&$\mathbf{1.09}$ &2.28&5.31&29.06&12.29&$--$&7.99\\
    ~&~&IT       &1403  &37 &1369&8075&1049&$--$&72\\ \cline{2-10}
    ~&\multirow{2}*{$k_{b}=50$}&CPU&$\mathbf{1.67}$ &5.63&14.03&72.82&31.04&$--$&19.72\\
    ~&~&IT       &1427  &37 &1363&81210&1047&$--$&73\\ \hline
    \multirow{6}*{randn(4000,600)}
   ~&\multirow{2}*{$k_{b}=10$}&CPU &$\mathbf{0.93}$ &1.44 &4.06&18.91&8.52&$--$&5.48\\
   ~&~&IT    &1786 &45&1745&10147&1375&$--$&94\\ \cline{2-10}
   ~&\multirow{2}*{$k_{b}=20$}&CPU&$\mathbf{1.51}$  &2.86 &8.28&36.86&16.91&$--$&10.48\\
   ~&~&IT       &1790  &45 &1731&10104&1375&$--$&94\\ \cline{2-10}
   ~&\multirow{2}*{$k_{b}=50$}&CPU&$\mathbf{2.28}$ &7.11 &20.87&92.61&43.06&$--$&26.06\\
    ~&~&IT       &1788  &45 &1743&10162&1388&$--$&93\\ \hline
    \multirow{6}*{randn(4000,700)}
   ~&\multirow{2}*{$k_{b}=10$}&CPU &$\mathbf{1.45}$  &1.84 &6.31&23.64&11.64&$--$&6.86\\
   ~&~&IT    &2188 &52 &2165&12389&1759&$--$&115\\ \cline{2-10}
    ~&\multirow{2}*{$k_{b}=20$}&CPU&$\mathbf{2.06}$  &3.68&12.37&46.53&23.44&$--$&13.58\\
    ~&~&IT       &2211 &52 & 2167&12422&2775&$--$&117\\ \cline{2-10}
    ~&\multirow{2}*{$k_{b}=50$}&CPU&$\mathbf{3.30}$ &9.25&31.56&116.10&59.03&$--$&34.53\\
   ~&~&IT       &2282 &53 &2171&12449&1781&$--$&116\\ \hline
\end{tabular}
\end{center}
\end{scriptsize}
\end{table}

$\mathbf{Example~ 4.4}$. In this example, the text matrices are from the University of Florida sparse matrix collection \cite{6}. Table 4.6 lists the details of these matrices. Similarly, the right-hand side $\mathbf{B}=\mathbf{AX}^{\star}\in R^{m\times k_{b}}$, where $\mathbf{X}^{\star}\in R^{n\times k_{b}}$ is generated by using the MATLAB function $randn(n,k_{b})$. Table 4.7 lists the numerical results obtained from the seven methods.

Again, we see from Table 4.7 that Algorithm 5 is superior to the other methods in terms of CPU time, while the RK, REABK and RBK methods do not work in many cases. Indeed, Algorithm 5 is about 2--5 times faster than GBK, and is about 10--20 times faster than SRK and GRK. Thus, we benefit from updating all the solutions simultaneously, and Algorithm 5 achieves much higher computational efficiency than the other randomized Kaczmarz methods for large-scale linear system with multiple right-hand sides.

\begin{table}[H]
 \begin{small}
 \begin{center}
 \caption{$\mathbf{Example ~4.4}$. Test matrices used in Example 4.1, where ``$(\cdot)^T$'' stands for the transpose of a matrix.}
   \vspace{0.2cm}
\begin{tabular}{|c|c|c|c|c}
\hline  Name&size $(m\times n)$&nnz&background\\
\hline $nemsemm1^{T}$&$75352\times 3945$&1053986&Linear Programming Problem\\
\hline $bibd\_16\_8^{T}$&$12870\times 120$&360360&Combinatorial Problem\\
\hline seymourl&$6316\times 4944$&38493&Linear Programming Problem\\
\hline  $gen4^{T}$&$6316\times 4298$&1537&Linear Programming Problem\\
\hline  $lp\_bnl2^{T}$&$4486\times 2324$&14996&Linear Programming Problem\\
\hline  $cq9^{T}$&$21534\times9278$&96653&Linear Programming Problem\\
\hline  $lp\_80bau3b^{T}$&$12061\times2262$&23264&Linear Programming Problem\\
\hline  $r05^{T}$&$9690\times5190$&104145&Linear Programming Problem\\
\hline  $rosen10^{T}$&$6152\times2056$&64192&Linear Programming Problem\\
\hline

\end{tabular}
\end{center}
\end{small}
\end{table}

\begin{table}
\begin{scriptsize}
\begin{center}
\caption
  {$\mathbf{Example ~4.4}$. Numerical results of the randomized block Kaczmarz-type methods for solving \eqref{1.1} on some large-sparse matrices, $tol=10^{-3}$, where the best ones in terms of CPU time are in bold. The sampling ratio is chosen as $\eta=0.01$ in Algorithm 5. Except for Algorithm 5, the value of ``IT" is rounded up by the average number of solving the $k_{b}$ linear systems with single right-hand side.}
\begin{tabular}{c|c|c|c|c|c|c|c|c|c}
	\hline
    Name&&&${\rm Algorithm~5}$&GBK&SRK&RK&GRK&REABK&RBK\\ \hline
     \multirow{6}*{$nemsemm1^{T}$}
    &\multirow{2}*{$k_{b}=10$}&CPU &$\mathbf{18.24}$ &33.83&126.11&$--$&241.47&$--$&$--$\\
    ~&~&IT    &4131&19 &4052&$--$&4372&$--$&$--$\\ \cline{2-10}
    ~&\multirow{2}*{$k_{b}=20$}&CPU&$\mathbf{23.42}$  &92.26&242.73&$--$&482.29&$--$&$--$\\
    ~&~&IT       &4130  &18&4046&$--$&4355&$--$&$--$\\ \cline{2-10}
    ~&\multirow{2}*{$k_{b}=50$}&CPU&$\mathbf{50.95}$  &208.72&607.97&$--$&1176.22&$--$&$--$\\
    ~&~&IT       &4136 &18 &4047&$--$&4360&$--$&$--$\\ \hline
    \multirow{6}*{bide\_16\_8}
    &\multirow{2}*{$k_{b}=10$}&CPU &$\mathbf{0.35}$ &0.97 &3.12&761.50&2.44&$--$&14.77\\
   ~&~&IT    &228 &7 &225&17156&196&$--$&373\\ \cline{2-10}
    ~&\multirow{2}*{$k_{b}=20$}&CPU&$\mathbf{0.41}$ &1.91 &5.93&1493.90&4.93&$--$&28.34\\
    ~&~&IT       &242  &7&225&1706&198&$--$&357\\ \cline{2-10}
    ~&\multirow{2}*{$k_{b}=50$}&CPU&$\mathbf{0.53 }$ &4.33 &14.84&3724.40&12.20&$--$&70.98\\
    ~&~&IT       &236 &7 &224&1702&196&$--$&358\\ \hline
    \multirow{6}*{seymourl}
    &\multirow{2}*{$k_{b}=10$}&CPU &$\mathbf{6.42}$&19.54 &30.24 &$--$&45.23&$--$&403.21\\
    ~&~&IT    &5167  &24 &5047&$--$&5321&$--$&1702\\ \cline{2-10}
    ~&\multirow{2}*{$k_{b}=20$}&CPU&$\mathbf{6.89}$ &39.33 &61.86&$--$&90.72&$--$&795.20\\
    ~&~&IT       &5075 &22 &4966&$--$&5331&$--$&1733\\ \cline{2-10}
    ~&\multirow{2}*{$k_{b}=50$}&CPU&$\mathbf{21.62}$ &94.97 &147.82&$--$&222.98&$--$&2064.70\\
    ~&~&IT       &5193 &22 &5021&$--$&5305&$--$&1803\\ \hline
     \multirow{6}*{lp\_bnl2$^{T}$}
    &\multirow{2}*{$k_{b}=10$}&CPU &$\mathbf{59.98}$&176.68 &206.68 &$--$&292.04&$--$&34534.00\\
    ~&~&IT    &11637 &786 &82952&$--$&57369&$--$&320433\\ \cline{2-10}
    ~&\multirow{2}*{$k_{b}=20$}&CPU&$\mathbf{88.88}$ &318.30 &460.44&$--$&572.47&$--$&$--$\\
    ~&~&IT       &129578 &754 &90310&$--$&57762&$--$&$--$\\ \cline{2-10}
    ~&\multirow{2}*{$k_{b}=50$}&CPU&$\mathbf{169.79}$ &790.27 &1081.81&$--$&1391.30&$--$&$--$\\
    ~&~&IT       &166582 &722 &86530&$--$&56034&$--$&$--$\\ \hline
    \multirow{6}*{cq9$^{T}$}
    &\multirow{2}*{$k_{b}=10$}&CPU &$\mathbf{95.04}$&239.71&435.68 &$--$&419.25&$--$&$--$\\
    ~&~&IT    &31876 &206 &38676&$--$&25417&$--$&$--$\\ \cline{2-10}
    ~&\multirow{2}*{$k_{b}=20$}&CPU&$\mathbf{320.53}$ &520.74 &1053.11&$--$&1345.90&$--$&$--$\\
    ~&~&IT       &41300 &222 &45852&$--$&615&$--$&$--$\\ \cline{2-10}
    ~&\multirow{2}*{$k_{b}=50$}&CPU&$\mathbf{648.91}$ &1222.10 &2536.60&$--$&2272.20&$--$&$--$\\
    ~&~&IT       &42668 &214&44428&$--$&27474&$--$&$--$\\ \hline
     \multirow{6}*{lp\_80bau3b$^{T}$}
    &\multirow{2}*{$k_{b}=10$}&CPU &$\mathbf{2.60}$&16.12&13.08 &$--$&21.74&$--$&1446.60\\
    ~&~&IT    &2136 &12 &2122&$--$&2409&$--$&12274\\ \cline{2-10}
    ~&\multirow{2}*{$k_{b}=20$}&CPU&$\mathbf{3.24}$ &32.28 &29.68&$--$&44.69&$--$&2708.80\\
    ~&~&IT       &2166 &13 &2136&$--$&2438&$--$&11581\\ \cline{2-10}
    ~&\multirow{2}*{$k_{b}=50$}&CPU&$\mathbf{4.94}$ &76.97&64.56&$--$&96.19&$--$&6565.40\\
    ~&~&IT       &2166 &13&2122&$--$&2463&$--$&11899\\ \hline
    \multirow{6}*{r05$^{T}$}
    &\multirow{2}*{$k_{b}=10$}&CPU &$\mathbf{57.07}$&212.93&247.95 &$--$&258.20&$--$&185.97\\
    ~&~&IT    &26940 &592 &24588&$--$&25005&$--$&859\\ \cline{2-10}
    ~&\multirow{2}*{$k_{b}=20$}&CPU&$\mathbf{68.11}$ &441.46 &573.79&$--$&507.72&$--$&5434.95\\
    ~&~&IT       &28655 &608&24719&$--$&24934&$--$&979\\ \cline{2-10}
    ~&\multirow{2}*{$k_{b}=50$}&CPU&$\mathbf{270.66}$ &1013.40&1241.10&$--$&1290.40&$--$&956.12\\
    ~&~&IT       &28814&611&24860&$--$&25174&$--$&882\\ \hline
     \multirow{6}*{rosen10$^{T}$}
    &\multirow{2}*{$k_{b}=10$}&CPU &$\mathbf{2.1}$&4.99&15.15 &$--$&16.04&$--$&18.50\\
    ~&~&IT    &2097 &41&2129&$--$&2595&$--$&177\\ \cline{2-10}
    ~&\multirow{2}*{$k_{b}=20$}&CPU&$\mathbf{3.40}$ &9.77 &30.00&$--$&31.67&$--$&39.15\\
    ~&~&IT       &2161 &41&2100&$--$&2618&$--$&186\\ \cline{2-10}
    ~&\multirow{2}*{$k_{b}=50$}&CPU&$\mathbf{5.12}$ &23.92&75.54&$--$&81.73&$--$&96.50\\
    ~&~&IT       &2204&40&2119&$--$&2618&$--$&182\\ \hline

\end{tabular}
\end{center}
\end{scriptsize}
\end{table}

\section{Conclusions}
Inspired by the semi-randomized Kaczmarz method with simple sampling \cite{12}, we propose two semi-randomized Kaczmarz methods with simple sampling for single right-hand and multiple right-hand sides, respectively. The key is to use a small part of the residual vector to construct the working rows for the semi-randomized block Kaczmarz method. The convergence of the proposed methods is discussed. Comprehensive numerical experiments demonstrate that our methods often outperform some popular randomized Kaczmarz methods and randomized block Kaczmarz methods for large-scale linear systems with single or multiple right-hand sides.


\section{Conflict of Interest Statement}
None of the authors have a conflict of interest to disclose.

\section{Acknowledgements}
This work is supported by the National Natural Science Foundation of China under grant 12271518, the Fujian Natural Science Foundation under grant 2023J01354, the Key Research and Development Project of Xuzhou Natural Science Foundation under grant KC22288, and the Open Project of Key Laboratory of Data Science and Intelligence Education of the Ministry of Education under grant DSIE202203.

\end{document}